# SPANNING TREE SIZE IN RANDOM BINARY SEARCH TREES

By Alois Panholzer and Helmut Prodinger

*Technische Universität Wien and University of the Witwatersrand*

This paper deals with the size of the spanning tree of $p$ randomly chosen nodes in a binary search tree. It is shown via generating functions methods, that for fixed $p$, the (normalized) spanning tree size converges in law to the Normal distribution. The special case $p=2$ reproves the recent result (obtained by the contraction method by Mahmoud and Neininger [*Ann. Appl. Probab.* **13** (2003) 253–276]), that the distribution of distances in random binary search trees has a Gaussian limit law. In the proof we use the fact that the spanning tree size is closely related to the number of passes in Multiple Quickselect. This parameter, in particular, its first two moments, was studied earlier by Panholzer and Prodinger [*Random Structures Algorithms* **13** (1998) 189–209]. Here we show also that this normalized parameter has for fixed $p$-order statistics a Gaussian limit law. For $p=1$ this gives the well-known result that the depth of a randomly selected node in a random binary search tree converges in law to the Normal distribution.

**1. Introduction.** In the papers [7] and [1] the distances between nodes in random search trees, respectively, random recursive trees were studied. It was proven in [7] that the (edge) distances $\Delta_n$ between two randomly selected nodes in random binary search trees of size $n$ are asymptotically normally (Gaussian) distributed, where the so-called random permutation model was used as the model of randomness for the trees. This means that every permutation of length $n$ is assumed to be equally likely when generating a binary search tree; furthermore, for selecting nodes, all $\binom{n}{2}$ pairs of nodes are assumed to be equally likely.

In [1] it was shown that the distribution of the distance $\Delta(i,n)$ between a fixed node $i$ and the node $n$ in a random recursive tree of size $n$ is (for a fixed ratio $\rho := \frac{i}{n}$ with $0<\rho<1$) asymptotically normally distributed.

A related parameter to the distance between two randomly selected nodes is the *Wiener index* of a graph, which is defined to be the sum of all distances between pairs of nodes in the considered graph. The Wiener index was studied for certain families of graphs and, although the scaled mean of this parameter must coincide with the mean distance of two randomly selected nodes, it turned out that the Wiener index was asymptotically *not* normally distributed for random recursive trees and random binary search trees (see [8] and [5]).

In this paper we concentrate on random binary search trees and study a natural extension of the distance between two randomly selected nodes, namely the size of the spanning tree of $p$ randomly chosen nodes in the tree. Again, we use the random permutation model for the generation of the binary search trees and also that all $\binom{n}{p}$ possibilities to select $p$ nodes in a tree of size $n$ are equally likely. The selection of the $p$ nodes will thus be independent of the chosen tree. The random variable $Y_{n,p}$, which counts the size of the spanning tree of $p$ randomly selected nodes in a random binary search tree of size $n$, is then a direct extension of $\Delta_n$, since the edge distance between two nodes is nothing else than the size of the spanning tree of these two nodes minus one and thus $\Delta_n \stackrel{\mathcal{L}}{=} Y_{n,2} - 1$, where $\stackrel{\mathcal{L}}{=}$ denotes equality in distribution.

In the mathematical analysis of $Y_{n,p}$ we use the fact that it is closely related to the random variable $X_{n,p}$, which counts the number of passes required in the Multiple Quickselect algorithm to find a random $p$-order statistic in a data file of length $n$ (see [9] and the references cited therein for a description of this divide and conquer algorithm); the natural probability model for the data is, that their ranks form a random permutation of $\{1, \dots, n\}$ and we assume further that all $\binom{n}{p}$ sets of $p$-order statistics $\{1 \le i_1 < \cdots < i_p \le n\}$ are equally likely. Then by well-known relations between binary search trees and Quicksort-like algorithms, $X_{n,p}$ is equal to the number of ascendants of $p$ randomly chosen nodes in a random binary search tree of size $n$ or, equivalently, to the size of the spanning tree, spanned by

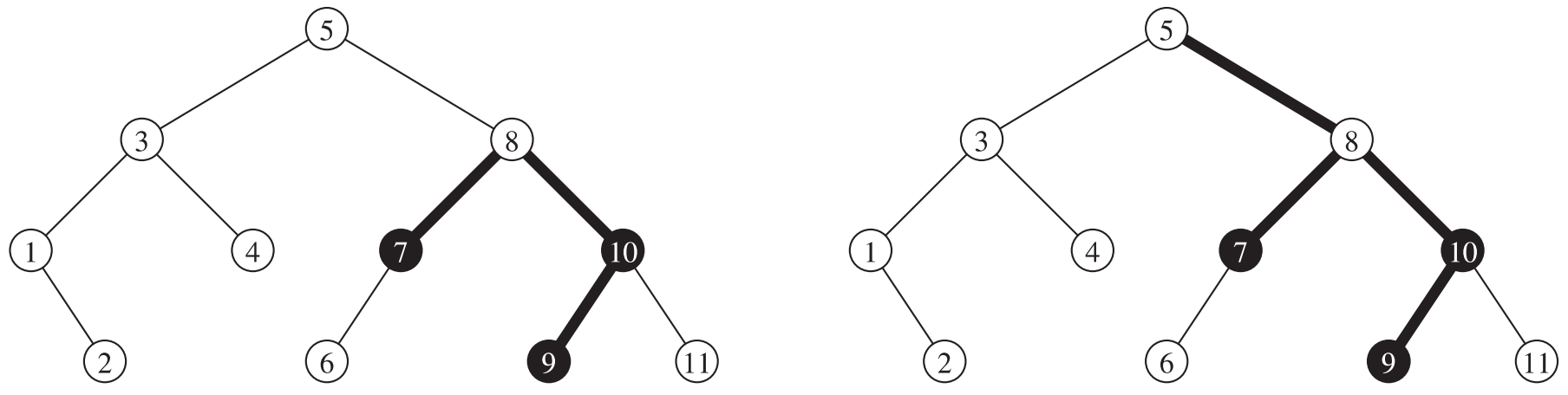

*Spanning tree of the nodes* 7, 9 *and* 10 *is of size* 4

5 *passes of the Multiple Quickselect algorithm to find the ranks* 7, 9 *and* 10

FIG. 1. *A binary search tree with the two parameters under consideration.*

the root and $p$ randomly chosen nodes (where of course the root could have been chosen as well) in a random binary search tree of size $n$. See Figure 1 for a comparison of both parameters.

The parameter $X_{n,p}$ was studied already in [9], where exact formulæ for the expectation and the variance were given. Here we show additionally that $X_{n,p}$ is, for fixed $p \geq 1$, asymptotically normally distributed (Philippe Flajolet mentioned that to Helmut Prodinger in 1998 without working out the details).

For $Y_{n,p}$ we also derive exact formulæ for the expectation and the variance and show that $Y_{n,p}$ is, for fixed $p \geq 2$, asymptotically normally distributed, where the special case $p = 2$ reproves that the distances $\Delta_n$ between two randomly selected nodes in random binary search trees of size $n$ are asymptotically normally distributed. Our approach uses generating functions, singularity analysis and a central limit theorem for combinatorial structures due to Hwang and avoids the difficulties which occur in [7] when showing the asymptotic normality for $\Delta_n$ using the contraction method which arises due to the degenerate nature of the distributional limit equation for $X_{n,1}$ (that was studied there to obtain the result for $\Delta_n$).

**2. Passes in Multiple Quickselect and spanning tree size in binary search trees.** First we want to translate the close relation between $X_{n,p}$ and $Y_{n,p}$ into an equation for suitable generating functions as described below.

Here we denote with $\varphi_{n,p,m} := \mathbb{P}\{X_{n,p} = m\}$ the probability that exactly $m$ passes of the Multiple Quickselect algorithm are required in order to find a random set of $p$-order statistics in a data file of length $n$ and with $F_{n,p,m} := \mathbb{P}\{Y_{n,p} = m\}$, the probability that the size of the spanning tree of $p$ randomly chosen nodes in a binary search tree of size $n$ is exactly $m$. Using the recursive structure of the search trees, we obtain for the generating functions $\phi_p(z,v) = \sum_{n,m\geq 0} \binom{n}{p} \varphi_{n,p,m} z^n v^m$, respectively, $F_p(z,v) = \sum_{n,m\geq 0} \binom{n}{p} F_{n,p,m} z^n v^m$ for $p \geq 1$ the recurrences

$$\frac{\partial}{\partial z}\phi_p(z,v) = v\sum_{i=0}^{p}\phi_i(z,v)\phi_{p-i}(z,v) + v\sum_{i=0}^{p-1}\phi_i(z,v)\phi_{p-1-i}(z,v) \tag{1}$$

and

$$\frac{\partial}{\partial z}F_p(z,v) = v\sum_{i=1}^{p-1}\phi_i(z,v)\phi_{p-i}(z,v) + v\sum_{i=0}^{p-1}\phi_i(z,v)\phi_{p-1-i}(z,v) + 2F_0(z,v)F_p(z,v), \tag{2}$$

with the initial functions $\phi_0(z,v) = F_0(z,v) = \frac{1}{1-z}$. The difference in the above recurrences reflects the difference between both parameters coming from the instance where the root is not selected and also the left $(i = 0)$, respectively, right ($i=p$) subtree of the root does not contain a selected node.

Introducing the trivariate generating functions $\Phi(z,u,v)=\sum_{p\ge 0}\phi_p(z,v)u^p$ and $F(z,u,v)=\sum_{p\ge 0}F_p(z,v)u^p$, we obtain first from (1) a Riccati differential equation

$$\frac{\partial}{\partial z}\Phi(z,u,v)=v(1+u)\Phi^2(z,u,v)+\frac{1-v}{(1-z)^2}, \tag{3}$$

with the initial value $\Phi(0,u,v)=1$. The solution of this equation is already given in [9],

$$\Phi(z,u,v)=\frac{\Omega+1-2v+(1-z)^{\Omega}(\Omega-1+2v)}{(\Omega+1-2v(1+u)+(1-z)^{\Omega}(\Omega-1+2v(1+u)))(1-z)}, \tag{4}$$

with

$$\Omega=\sqrt{1-4(1+u)v(1-v)}. \tag{5}$$

For $F(z,u,v)$ we get from (2) the differential equation

$$\frac{\partial}{\partial z}F(z,u,v) = v(1+u)\Phi^2(z,u,v)-\frac{2v}{1-z}\Phi(z,u,v)+\frac{2}{1-z}F(z,u,v)+\frac{v-1}{(1-z)^2}$$

or

$$\frac{\partial}{\partial z}F(z,u,v) = \frac{\partial}{\partial z}\Phi(z,u,v)-\frac{2v}{1-z}\Phi(z,u,v)+\frac{2}{1-z}F(z,u,v)+\frac{2(v-1)}{(1-z)^2},$$

with $F(0,u,v)=1$.

This equation then has the solution

$$\begin{aligned} F(z,u,v) &= \frac{1+2z(v-1)}{(1-z)^2} \\ &\quad + \frac{1}{(1-z)^2}\int_0^z\left[\frac{\partial}{\partial t}\Phi(t,u,v)-\frac{2v}{1-t}\Phi(t,u,v)\right](1-t)^2\,dt, \end{aligned} \tag{6}$$

with $\Phi(z,u,v)$ given by (4).

**3. Expectation and variance of the spanning tree size.** From (6) it is easy to obtain exact formulæ for the expectation

$$E_{n,p}=\mathbb{E}(Y_{n,p})=\frac{1}{\binom{n}{p}}[z^nu^p]\frac{\partial}{\partial v}F(z,u,v)\Big|_{v=1}$$

and the second factorial moment

$$M_{n,p}^{(2)} = \mathbb{E}(Y_{n,p}(Y_{n,p}-1)) = \frac{1}{\binom{n}{p}}[z^n u^p] \frac{\partial^2}{\partial v^2} F(z,u,v)\Big|_{v=1}$$

(and thus also for the variance $V_{n,p}$) of $Y_{n,p}$, the size of the spanning tree of $p$ randomly selected nodes in a random binary search tree of size $n$.

Differentiating (6) with respect to $v$ and evaluating at $v=1$ gives the following equation for $E(z,u) := \frac{\partial}{\partial v} F(z,u,v)|_{v=1}$:

$$\begin{aligned} E(z,u) = \frac{2z}{(1-z)^2} + \frac{1}{(1-z)^2} \int_0^z \Bigg[ & \frac{4(1-t)(1+u)u^2}{(1-t(1+u))^3} \log \frac{1}{1-t} \\ & + \frac{X}{(1-t(1+u))^3} \Bigg] dt, \end{aligned} \tag{7}$$

with $X = (-2+u) + (6+3u-3u^2-4u^3)t + (1+u)(2u^2-3u-6)t^2 + (2+u)(1+u)^2 t^3$.

This can be simplified to

$$\begin{aligned} E(z,u) = & \frac{2u(1+u)}{(1-z(1+u))^2} \log \frac{1}{1-z} - \frac{2u}{(1-z)^2(1+u)^2} \log \frac{1}{1-z(1+u)} \\ & + \frac{zu(1-2z-3u+z^2+uz-2u^2+2uz^2+3u^2z+u^2z^2)}{(1-z)^2(1+u)(1-z(1+u))^2}. \end{aligned} \tag{8}$$

To extract coefficients we use here and in the sequel the general formulæ (see, e.g., [3])

$$[z^n] \frac{1}{(1-z)^{m+1}} \log \frac{1}{1-z} = \binom{n+m}{n} (H_{n+m} - H_m),$$

$$[z^n] \frac{1}{(1-z)^{m+1}} \log^2 \frac{1}{1-z} = \binom{n+m}{n} ((H_{n+m} - H_m)^2 - (H_{n+m}^{(2)} - H_m^{(2)})),$$

where $H_n = \sum_{k=1}^n \frac{1}{k}$ and $H_n^{(2)} = \sum_{k=1}^n \frac{1}{k^2}$ denote the first and second order harmonic numbers.

By lengthy, but routine calculations, we finally get for $E_{n,p} = \frac{1}{\binom{n}{p}}[z^n u^p] E(z,u)$ an exact formula, which is given in the next lemma:

Lemma 1. *The expectation $E_{n,p} = \mathbb{E}(Y_{n,p})$ of the size of the spanning tree of $p$ randomly chosen nodes in a random binary search tree of size $n$ is for $p \geq 1$ given by*

$$\begin{aligned} E_{n,p} = & \frac{2p(n+1)^2}{(n+2-p)(n+1-p)}(H_n - H_p) + \frac{2(2p-1)(n+1)}{(n+2-p)(n+1-p)} + 3 + 2p \\ & - \frac{2pn}{n+1-p} + \frac{2p(n+1)(-1)^p}{\binom{n}{p}} H_n + \frac{2p(n+1)(-1)^p}{\binom{n}{p}} \sum_{k=1}^{p-1} \frac{(-1)^k}{k} \binom{n}{k}, \end{aligned}$$

*and asymptotically for fixed $p \geq 2$ by*

$$E_{n,p} = 2p \log n + 2p\gamma - 2pH_p + 3 - 2p - \frac{2p}{p-1} + \mathcal{O}\Big(\frac{\log n}{n}\Big).$$

For $p=1$, the formula simplifies to $E_{n,1}=1$ as it should.

We remark that

$$H_n = \sum_{k=1}^{n} \frac{(-1)^{k-1}}{k} \binom{n}{k},$$

and so one can give the alternative formula

$$E_{n,p} = \frac{2p(n+1)^2}{(n+2-p)(n+1-p)}(H_n - H_p) + \frac{2(2p-1)(n+1)}{(n+2-p)(n+1-p)} + 3 + 2p - \frac{2pn}{n+1-p} + \frac{2p(n+1)(-1)^p}{\binom{n}{p}} \sum_{k=p}^{n} \frac{(-1)^{k-1}}{k}\binom{n}{k}.$$

When we differentiate equation (6) twice with respect to $v$ and evaluate at $v=1$, we finally obtain the following formula for $M_2(z,u) := \frac{\partial^2}{\partial v^2} F(z,u,v)|_{v=1}$:

$$\begin{aligned} M_2(z,u) = &-\frac{8u}{(1-z)^2}\int_0^z \frac{1}{1-t(1+u)} \log\frac{1}{1-t}\,dt \\ &+ \frac{4u(1+u)^2(1-z+2u-uz)}{(1-z(1+u))^3}\log^2\frac{1}{1-z} \\ &+ \frac{12u}{(1-z)^2(1+u)^2}\log\frac{1}{1-z(1+u)} \\ &+ \frac{4u\Psi_1}{(1-z)(1-z(1+u))^3}\log\frac{1}{1-z} \\ &+ \frac{2u^2 z\Psi_2}{(1-z)^2(1-z(1+u))^3(1+u)}, \end{aligned} \tag{9}$$

with the abbreviations

$$\begin{aligned} \Psi_1 &= -z^2u^2 - 3z^2 - 5z^2u + u^3z^2 + 6z + 9zu + 2u^2z - u^3z - 3u^2 - 4u - 3, \\ \Psi_2 &= -z^3u^3 + 3u^3z^2 - 2u^3z - 19u^2z + 22z^2u^2 + 6u^2 - 3z^3u^2 \\ &\quad + 14u - 3z^3u - 46zu + 35z^2u + 14 - 29z - z^3 + 16z^2. \end{aligned}$$

Extracting coefficients gives after a somewhat lengthy calculation an exact formula for the second factorial moment $M_{n,p}^{(2)} = \frac{1}{\binom{n}{p}}[z^n u^p]M_2(z,u)$ and we get via $V_{n,p} = M_{n,p}^{(2)} + E_{n,p} - E_{n,p}^2$ the following result:

LEMMA 2. *The variance $V_{n,p} = \mathbb{V}(Y_{n,p})$ of the size of the spanning tree of $p$ chosen nodes in a random binary search tree of size $n$ is for $p \geq 2$ given by*

$$\begin{aligned}
V_{n,p} &= \frac{4(-1)^p(n+1)(2pH_n - 2pH_p + 2 - 3p^2)}{p\binom{n}{p}} \sum_{k=1}^{p-1} \frac{(-1)^k}{k}\binom{n}{k} \\
&- \frac{8(-1)^p(n+1)}{\binom{n}{p}} \sum_{k=1}^{p-1} \frac{(-1)^k}{k^2}\binom{n}{k} \\
&+ \frac{4(-1)^p(n+1)}{\binom{n}{p}}(H_n^2 - H_n^{(2)} - 2H_pH_n) \\
&+ \frac{4(-1)^p(2-3p^2)(n+1)}{p\binom{n}{p}} H_n - \frac{4\Psi_3}{(n+4-p)^{\underline{4}}}(H_n - H_p) \\
&+ \frac{4p(n+2)(n+1)^3(np+2+p)}{(n+4-p)^{\underline{4}}}((H_n - H_p)^2 - (H_n^{(2)} - H_p^{(2)})) \\
&+ \frac{2\Psi_4}{(n+4-p)^{\underline{4}}} + E_{n,p} - E_{n,p}^2,
\end{aligned}$$

*with*

$$\begin{aligned}
\Psi_3 &= -2p^4 - 6n^2p^3 + 16p^3 - 2n^3p^3 - 45p^2n - 58p^2 - 4n^2p^2 + 2p^2n^4 \\
&\quad + 7n^3p^2 + 56p + 78np + 6n^3p + 41pn^2 - n^4p - 8 - 20n - 16n^2 - 4n^3, \\
\Psi_4 &= -144 - 6p^5 - 3n^4 - 152p^3n + 2p^4n^2 + 25p^4n - 234n + 78p + 10np \\
&\quad - 5n^4p - 39n^2p^3 - 4n^3p^3 + 250p^2n + 119n^2p^2 + 2p^2n^4 + 25n^3p^2 \\
&\quad - 22n^3p - 35pn^2 + 155p^2 - 173p^3 + 58p^4 - 153n^2 - 42n^3.
\end{aligned}$$

*Further we have $V_{n,1} = 0$ and the following asymptotic expansion for $n \to \infty$ and fixed $p \geq 2$:*

$$\begin{aligned}
V_{n,p} &= 2p\log n - 2p(H_p - \gamma) - 4p^2\left(\frac{\pi^2}{6} - H_p^{(2)}\right) \\
&\quad + \frac{2(-2+7p-5p^2+2p^3)}{(1-p)^2} + \mathcal{O}\left(\frac{\log^2 n}{n}\right).
\end{aligned}$$

Here we used the abbreviation $x^{\underline{m}} := x(x-1)\cdots(x-m+1)$ for the falling factorials.

We remark again that an alternative representation of the variance would be possible using the additional formula

$$\frac{1}{2}(H_n^2 + H_n^{(2)}) = \sum_{k=1}^{n} \frac{(-1)^{k-1}}{k^2}\binom{n}{k}.$$

**4. The limiting distribution of the number of passes in Multiple Quickselect.** We will show that both random variables $X_{n,p}$ and $Y_{n,p}$ satisfy, for fixed $p$, a Gaussian limit law. To do this, we will expand the coefficients at $u^p$ (for fixed $p$) of the trivariate generating functions $\Phi(z,u,v)$, respectively, $F(z,u,v)$ around their dominant singularity $z=1$, where the expansion holds uniformly for $|v-1|<\tau$, for $\tau>0$. Singularity analysis (see [2]) of generating functions allows then to translate these expansions into an asymptotic expansion of the moment generating function (the Laplace transform) of the considered random variables. Then we can apply the so called Quasi power theorem (see [4]) to establish the weak convergence of the random variables to the normal distribution with certain convergence rates.

In this section we will treat the random variable $X_{n,p}$. As described above, we are interested in an asymptotic expansion of $\frac{1}{\binom{n}{p}}[z^n u^p]\Phi(z,u,v)$ for $n\to\infty$ and fixed $p$ uniformly for $|v-1|\le\tau$, where $\Phi(z,u,v)$ is given by the exact formula (4).

To expand $\Phi(z,u,v)$ we will use some auxiliary expansions of

$$f(u) = \Omega + 1 - 2v + (1-z)^{\Omega}(\Omega - 1 + 2v), \tag{10}$$

$$g(u) = \Omega + 1 - 2v(1+u) + (1-z)^{\Omega}(\Omega - 1 + 2v(1+u)) \tag{11}$$

with $\Omega$ given by (5). All $\mathcal{O}$-terms in the expansions given below are uniform for $|v-1|\le\tau$, as required. In the sequel, we will use the notations $D_u$ for the differential operator w.r.t. $u$ and $N_u$ for the evaluation operator at $u=0$.

Since

$$\Omega = (2v-1)\sum_{k\ge 0}\binom{1/2}{k}\left(-\frac{4v(1-v)}{1-4v(1-v)}\right)^k u^k,$$

we get

$$\begin{aligned}(1-z)^{\Omega} &= e^{\Omega\log(1-z)}\\ &= (1-z)^{2v-1}\\ &\quad\times\exp\left[(2v-1)\log(1-z)\sum_{k\ge 1}\binom{1/2}{k}\left(-\frac{4v(1-v)}{1-4v(1-v)}\right)^k u^k\right],\end{aligned}$$

and thus

$$N_u D_u^p (1-z)^{\Omega} = \mathcal{O}((1-z)^{2v-1}\log^p(1-z)). \tag{12}$$

We have

$$f(0) = g(0) = (1-z)^{2v-1}2(2v-1) \tag{13}$$

and we get further

$$N_u D_u f(u) = -\frac{4v(1-v)}{2(2v-1)} + \mathcal{O}((1-z)^{2v-1}\log(1-z)), \tag{14a}$$

$$N_u D_u g(u) = -\frac{2v^2}{2v-1} + \mathcal{O}((1-z)^{2v-1}\log(1-z)), \tag{14b}$$

$$\begin{aligned} N_u D_u^p f(u) &= (2v-1)p!\binom{1/2}{p}\left(-\frac{4v(1-v)}{1-4v(1-v)}\right)^p \\ &\quad + \mathcal{O}((1-z)^{2v-1}\log^p(1-z)), \end{aligned} \tag{14c}$$

$$\begin{aligned} N_u D_u^p g(u) &= (2v-1)p!\binom{1/2}{p}\left(-\frac{4v(1-v)}{1-4v(1-v)}\right)^p \\ &\quad + \mathcal{O}((1-z)^{2v-1}\log^p(1-z)), \end{aligned} \tag{14d}$$

for $p \geq 2$. Furthermore, we want to expand $D_u^p(g(u))^{-1}$ (for $p \geq 1$) in terms of falling powers of $(g(u))^{-1}$, which gives

$$\begin{aligned} &D_u^p(g(u))^{-1} \\ &\qquad = (-1)^p p!(g(u))^{-p-1}(g'(u))^p \\ &\qquad\quad + \frac{(-1)^{p-1}(p-1)p!}{2}(g(u))^{-p}(g'(u))^{p-2}g''(u) + \mathcal{O}((g(u))^{-p+1}) \end{aligned}$$

and hence we obtain the expansion

$$\begin{aligned} &N_u D_u^p(g(u))^{-1} \\ &\qquad = (-1)^p p!\frac{(-2v^2/(2v-1) + \mathcal{O}((1-z)^{2v-1}\log(1-z)))^p}{(1-z)^{(p+1)(2v-1)}2^{p+1}(2v-1)^{p+1}} \\ &\qquad\quad + \mathcal{O}\left(\frac{1}{(1-z)^{p(2v-1)}}\right) \\ &\qquad = \frac{(-1)^p p!(-2v^2/(2v-1))^p}{(1-z)^{(p+1)(2v-1)}2^{p+1}(2v-1)^{p+1}} + \mathcal{O}\left(\frac{\log(1-z)}{(1-z)^{p(2v-1)}}\right). \end{aligned} \tag{15}$$

This finally gives

$$\begin{aligned} &N_u D_u^p \Phi(z,u,v) \\ &\qquad = N_u D_u^p \frac{f(u)}{(1-z)g(u)} \\ &\qquad = \frac{1}{1-z}f(0)N_u D_u^p(g(u))^{-1} + \frac{1}{1-z}pf'(0)N_u D_u^{p-1}(g(u))^{-1} \\ &\qquad\quad + \frac{1}{1-z}\mathcal{O}\left(\frac{1}{(1-z)^{(p-1)(2v-1)}}\right) \end{aligned}$$

$$
\begin{aligned}
&= \frac{(1-z)^{2v-1}2(2v-1)(-1)^p p!(-2v^2/(2v-1))^p}{(1-z)^{(p+1)(2v-1)+1}2^{p+1}(2v-1)^{p+1}} \\
&\quad + \frac{p(-4v(1-v)/(2(2v-1)))(-1)^{p-1}(p-1)!(-2v^2/(2v-1))^{p-1}}{(1-z)^{p(2v-1)+1}2^p(2v-1)^p} \\
&\quad + \mathcal{O}\Big(\frac{\log(1-z)}{(1-z)^{(p-1)(2v-1)+1}}\Big) \\
&= \frac{p!(v/(2v-1))^{2p-1}}{(1-z)^{p(2v-1)+1}} + \mathcal{O}\Big(\frac{\log(1-z)}{(1-z)^{(p-1)(2v-1)+1}}\Big).
\end{aligned} \tag{16}
$$

Singularity analysis leads then directly to

$$
\begin{aligned}
&[z^n]N_u D_u^p \Phi(z,u,v) \\
&\quad = p!\Big(\frac{v}{2v-1}\Big)^{2p-1}\frac{n^{p(2v-1)}}{\Gamma(p(2v-1)+1)}\Big(1+\mathcal{O}\Big(\frac{1}{n}\Big)\Big) \\
&\qquad + \mathcal{O}((\log n) n^{(p-1)(2v-1)}) \\
&\quad = \frac{p!(v/(2v-1))^{2p-1}n^{p(2v-1)}}{\Gamma(p(2v-1)+1)}\Big(1+\mathcal{O}\Big(\frac{1}{n}\Big)\Big)\Big(1+\mathcal{O}\Big(\frac{\log n}{n^{2v-1}}\Big)\Big) \\
&\quad = \frac{p!(v/(2v-1))^{2p-1}n^{p(2v-1)}}{\Gamma(p(2v-1)+1)}\Big(1+\mathcal{O}\Big(\frac{1}{n^{1-\varepsilon}}\Big)\Big),
\end{aligned} \tag{17}
$$

uniformly for $\varepsilon > 0$ and $|v-1| \le \tau := \frac{\varepsilon}{3}$ and also to the following expansion, which is valid for fixed $p \ge 1$:

$$
\begin{aligned}
&\frac{1}{\binom{n}{p}}[z^n u^p]\Phi(z,u,v) \\
&\quad = \frac{p!}{n^p}[z^n u^p]\Phi(z,u,v)\Big(1+\mathcal{O}\Big(\frac{1}{n}\Big)\Big) \\
&\quad = \frac{1}{n^p}[z^n]N_u D_u^p \Phi(z,u,v)\Big(1+\mathcal{O}\Big(\frac{1}{n}\Big)\Big) \\
&\quad = \frac{p!(v/(2v-1))^{2p-1}n^{p(2v-2)}}{\Gamma(p(2v-1)+1)}\Big(1+\mathcal{O}\Big(\frac{1}{n^{1-\varepsilon}}\Big)\Big) \\
&\quad = \exp\Big[p(2v-2)\log n + \log\Big(\frac{p!(v/(2v-1))^{2p-1}}{\Gamma(p(2v-1)+1)}\Big)\Big] \\
&\qquad \times \Big(1+\mathcal{O}\Big(\frac{1}{n^{1-\varepsilon}}\Big)\Big).
\end{aligned} \tag{18}
$$

We give here the Quasi power theorem as proven in [4], which we want to apply to our problem.

THEOREM 3 (H. K. Hwang). *Let $\{\Omega_n\}_{n\geq 1}$ be a sequence of integral random variables. Suppose that the moment generating function satisfies the asymptotic expression*

$$M_n(s) := \mathbb{E}(e^{\Omega_n s}) = \sum_{m\geq 0} \mathbb{P}\{\Omega_n = m\} e^{ms} = e^{H_n(s)}(1 + \mathcal{O}(\kappa_n^{-1})),$$

*the $\mathcal{O}$-term being uniform for $|s| \leq \tau$, $s \in \mathbb{C}$, $\tau > 0$, where:*

(i) $H_n(s) = u(s)\phi(n) + v(s)$, *with $u(s)$ and $v(s)$ analytic for $|s| \leq \tau$ and independent of $n$; $u''(0) \neq 0$,*
(ii) $\phi(n) \to \infty$,
(iii) $\kappa_n \to \infty$.

*Under these assumptions, the distribution of $\Omega_n$ is asymptotically Gaussian*

$$\mathbb{P}\left\{\frac{\Omega_n - u'(0)\phi(n)}{\sqrt{u''(0)\phi(n)}} < x\right\} = \Phi(x) + \mathcal{O}\left(\frac{1}{\kappa_n} + \frac{1}{\sqrt{\phi(n)}}\right),$$

*uniformly with respect to $x$, $x \in \mathbb{R}$. Here $\Phi(x)$ denotes the distribution function of the standard normal distribution $\mathcal{N}(0,1)$.*

*Moreover, the mean and the variance of $\Omega_n$ satisfy*

$$\mathbb{E}(\Omega_n) = u'(0)\phi(n) + v'(0) + \mathcal{O}(\kappa_n^{-1}),$$
$$\mathbb{V}(\Omega_n) = u''(0)\phi(n) + v''(0) + \mathcal{O}(\kappa_n^{-1}).$$

From (18) we get, with the notation in Theorem 3,

$$u(s) = p(2e^s - 2), \qquad v(s) = \log\left(\frac{p!(e^s/(2e^s-1))^{2p-1}}{\Gamma(p(2e^s-1)+1)}\right),$$
$$\phi(n) = \log n, \qquad \kappa_n = n^{1-\varepsilon}.$$

We find

$$u'(0) = 2p, \qquad u''(0) = 2p, \tag{19}$$

and

$$\begin{aligned} v'(0) &= -2p + 1 - 2p\Psi(p+1) = -2pH_p + 2p\gamma + 1 - 2p, \\ v''(0) &= 2(2p-1) - 2p\Psi(p+1) - 4p^2\Psi'(p+1) \\ &= 2(2p-1) - 2pH_p + 2p\gamma - \tfrac{2}{3}\pi^2p^2 + 4p^2H_p^{(2)}, \end{aligned} \tag{20}$$

where $\Psi(x)$ denotes the digamma function: $\Psi(x) := (\log\Gamma(x))'$.

Hence, with equations (19) and (20), we get from Theorem 3 the following result:

THEOREM 4. *The distribution of the random variable $X_{n,p}$, which counts the number of passes in the Multiple Quickselect algorithm that are required to find a random order statistic of $p$ elements in a data file of size $n$, is for fixed $p \ge 1$ asymptotically Gaussian, where the convergence rate is of order $\mathcal{O}(\frac{1}{\sqrt{\log n}})$:*

$$\mathbb{P}\left\{\frac{X_{n,p} - 2p\log n}{\sqrt{2p\log n}} < x\right\} = \Phi(x) + \mathcal{O}\left(\frac{1}{\sqrt{\log n}}\right),$$

*and the expectation $E_{n,p} = \mathbb{E}(X_{n,p})$ and the variance $V_{n,p} = \mathbb{V}(X_{n,p})$ satisfy*

$$E_{n,p} = 2p\log n + 1 - 2p - 2pH_p + 2p\gamma + \mathcal{O}\left(\frac{1}{n^{1-\varepsilon}}\right),$$

$$V_{n,p} = 2p\log n + 4p - 2 - 2pH_p + 2p\gamma + 4p^2H_p^{(2)} - \frac{2}{3}\pi^2p^2 + \mathcal{O}\left(\frac{1}{n^{1-\varepsilon}}\right).$$

The result for $E_{n,p}$ and $V_{n,p}$ already appeared in [9], but unfortunately there was a typo in the formula for $V_{n,p}$.

For the case $p = 1$ we have that $X_{n,1}$ counts the number of comparisons encountered by a successful search in a random binary search tree and this is, up to an additive constant, the same as the depth $D_n$ of a randomly selected node, $D_n \stackrel{\mathcal{L}}{=} X_{n,1} - 1$. The asymptotic normality of the distribution of $X_{n,1}$ is well known (see, e.g., [6]) and the convergence rate was recently established in [7].

**5. The limiting distribution of the spanning tree size in binary search trees.** In this section we will show that the normalized random variable $Y_{n,p}$, as defined in Section 4, has for fixed $p$ a Gaussian limiting distribution. Hence we are interested in an asymptotic expansion of $\frac{1}{\binom{n}{p}}[z^nu^p]F(z,u,v)$ for $n \to \infty$ and fixed $p$ uniformly for $|v-1| \le \tau$, where $F(z,u,v)$ is given by (6).

To do this, we will first study the behavior near the singularity $z = 1$ of the expression

$$\tilde{\Phi}(z,u,v) = \left(\frac{\partial}{\partial z}\Phi(z,u,v) - \frac{2v}{1-z}\Phi(z,u,v)\right)(1-z)^2, \tag{21}$$

which we can write as

$$\tilde{\Phi}(z,u,v) = \frac{\tilde{f}(u)}{(g(u))^2}, \tag{22}$$

where the function $\tilde{f}(u)$ is defined by

$$\begin{aligned}\tilde{f}(u) &= -\Omega(\Omega - 1 + 2v)(1-z)^{\Omega}g(u) \\ &\quad + \Omega(\Omega - 1 + 2v(1+u))(1-z)^{\Omega}f(u) + (1-2v)f(u)g(u),\end{aligned} \tag{23}$$

and $\Phi(z,u,v)$, $\Omega$, $f(u)$ and $g(u)$ are given by equations (4), (5), (10) and (11), respectively.

The relevant expansions are now

$$\tilde{f}(0) = -4(2v-1)^3(1-z)^{4v-2},$$
$$\tilde{f}'(0) = 8v^2(2v-1)(1-z)^{2v-1} + \mathcal{O}(\log(1-z)(1-z)^{4v-2}),$$
$$\tilde{f}''(0) = \frac{8(v-1)v^3}{2v-1} + \mathcal{O}(\log(1-z)(1-z)^{2v-1})$$

and

$$N_u D_u^p (g(u))^{-2} = (-1)^p (p+1)! \frac{(g'(0))^p}{(g(0))^{p+2}} + \mathcal{O}\left(\frac{1}{(g(0))^{p+1}}\right),$$

which leads, for $p \geq 2$, eventually to

$$(24)\qquad \begin{aligned} N_u D_u^p \tilde{\Phi}(z,u,v) &= \tilde{f}(0) N_u D_u^p (g(u))^{-2} + p\tilde{f}'(0) N_u D_u^{p-1}(g(u))^{-2} \\ &\quad + \frac{p(p-1)}{2}\tilde{f}''(0) N_u D_u^{p-2}(g(u))^{-2} + \mathcal{O}\left(\frac{1}{(1-z)^{(p-1)(2v-1)}}\right) \\ &= \frac{(p-1)p!v(v/(2v-1))^{2p-2}}{(1-z)^{p(2v-1)}} + \mathcal{O}\left(\frac{\log(1-z)}{(1-z)^{(p-1)(2v-1)}}\right). \end{aligned}$$

This gives then

$$(25)\qquad \begin{aligned} \frac{1}{\binom{n}{p}}[z^n u^p] F(z,u,v) &= \frac{1}{\binom{n}{p}}[z^n u^p] \frac{1}{(1-z)^2} \int_{t=0}^{z} \tilde{\Phi}(t,u,v)\, dt \\ &= \frac{1}{n^p}[z^n] \frac{1}{(1-z)^2} \int_{t=0}^{z} N_u D_u^p \tilde{\Phi}(t,u,v)\, dt \left(1 + \mathcal{O}\left(\frac{1}{n}\right)\right) \\ &= \frac{1}{n^p}[z^n] \frac{1}{(1-z)^2} \\ &\quad \times \int_{t=0}^{z} \left[\frac{(p-1)p!v(v/(2v-1))^{2p-2}}{(1-z)^{p(2v-1)}} \right. \\ &\qquad \left. + \mathcal{O}\left(\frac{\log(1-z)}{(1-z)^{(p-1)(2v-1)}}\right)\right] dt \left(1 + \mathcal{O}\left(\frac{1}{n}\right)\right). \end{aligned}$$

We get via singularity analysis

$$[z^n] \frac{1}{(1-z)^2} \int_{t=0}^{z} \frac{(p-1)p!v(v/(2v-1))^{2p-2}}{(1-z)^{p(2v-1)}}\, dt$$

$$= [z^n]\frac{(p-1)p!v(v/(2v-1))^{2p-2}}{(p(2v-1)-1)(1-z)^{p(2v-1)+1}} \tag{26}$$

$$= \frac{(p-1)p!v(v/(2v-1))^{2p-2}n^{p(2v-1)}}{(p(2v-1)-1)\Gamma(p(2v-1)+1)}\left(1+\mathcal{O}\left(\frac{1}{n}\right)\right)$$

and

$$[z^n]\frac{1}{(1-z)^2}\int_{t=0}^{z}\mathcal{O}\left(\frac{\log(1-t)}{(1-t)^{(p-1)(2v-1)}}\right)dt$$

$$= \mathcal{O}\left(\sum_{k=1}^{n}[z^{n-k}]\frac{1}{(1-z)^2}[z^k]\int_{t=0}^{z}\frac{\log(1-t)}{(1-t)^{(p-1)(2v-1)}}\,dt\right)$$

$$= \mathcal{O}\left(n\max_{1\le k\le n}[z^{n-k}]\frac{1}{(1-z)^2}\max_{1\le k\le n}[z^k]\int_{t=0}^{z}\frac{\log(1-t)}{(1-t)^{(p-1)(2v-1)}}\,dt\right) \tag{27}$$

$$= \mathcal{O}\left(n^2[z^n]\int_{t=0}^{z}\frac{\log(1-t)}{(1-t)^{(p-1)(2v-1)}}\,dt\right) = \mathcal{O}((\log n)\,n^{(p-1)(2v-1)})$$

$$= \mathcal{O}\left(n^{p(2v-1)}\frac{1}{n^{1-\varepsilon}}\right),$$

uniformly for $\varepsilon > 0$ and $|v-1| \le \tau := \frac{\varepsilon}{3}$.

Thus we obtain by combining the results (25)–(27) for $p \ge 2$ the asymptotic expansion

$$\frac{1}{\binom{n}{p}}[z^n u^p]F(z,u,v)$$

$$= \frac{(p-1)p!v(v/(2v-1))^{2p-2}n^{p(2v-2)}}{(p(2v-1)-1)\Gamma(p(2v-1)+1)}\left(1+\mathcal{O}\left(\frac{1}{n^{1-\varepsilon}}\right)\right). \tag{28}$$

To apply the Quasi power theorem, we write (28) as

$$\frac{1}{\binom{n}{p}}[z^n u^p]F(z,u,v)$$

$$= \exp\left[p(2v-2)\log n + \log\left(\frac{(p-1)p!v(v/(2v-1))^{2p-2}}{(p(2v-1)-1)\Gamma(p(2v-1)+1)}\right)\right] \tag{29}$$

$$\times\left(1+\mathcal{O}\left(\frac{1}{n^{1-\varepsilon}}\right)\right)$$

and then get, with the notation used in Theorem 3,

$$u(s) = p(2e^s-2), \qquad v(s) = \log\left(\frac{(p-1)p!e^s(e^s/(2e^s-1))^{2p-2}}{(p(2e^s-1)-1)\Gamma(p(2e^s-1)+1)}\right),$$

$$\phi(n) = \log n, \qquad \kappa_n = n^{1-\varepsilon}.$$

We have

$$u'(0) = 2p, \qquad u''(0) = 2p, \tag{30}$$

and

$$\begin{aligned} v'(0) &= -2p\Psi(p+1) + 3 - 2p - \frac{2p}{p-1} \\ &= -2pH_p + 2p\gamma + 3 - 2p - \frac{2p}{p-1}, \\ v''(0) &= -2p\Psi(p+1) - 4p^2\Psi'(p+1) + \frac{2(2p^3 - 5p^2 + 7p - 2)}{(p-1)^2} \\ &= -2pH_p + 2p\gamma - \frac{2}{3}\pi^2 p^2 + 4p^2 H_p^{(2)} + \frac{2(2p^3 - 5p^2 + 7p - 2)}{(p-1)^2}, \end{aligned} \tag{31}$$

which leads now to the following result:

THEOREM 5. *The distribution of the random variable $Y_{n,p}$, which counts the size of the spanning tree of $p$ randomly chosen nodes in a binary search tree of size $n$, is for fixed $p \geq 2$ asymptotically Gaussian, where the convergence rate is of order $\mathcal{O}(\frac{1}{\sqrt{\log n}})$.*

$$\mathbb{P}\left\{\frac{Y_{n,p} - 2p\log n}{\sqrt{2p\log n}} < x\right\} = \Phi(x) + \mathcal{O}\left(\frac{1}{\sqrt{\log n}}\right)$$

*and the expectation $E_{n,p} = \mathbb{E}(Y_{n,p})$ and the variance $V_{n,p} = \mathbb{V}(Y_{n,p})$ satisfy*

$$\begin{aligned} E_{n,p} &= 2p\log n - 2pH_p + 2p\gamma + 3 - 2p - \frac{2p}{p-1} + \mathcal{O}\left(\frac{1}{n^{1-\varepsilon}}\right), \\ V_{n,p} &= 2p\log n - 2pH_p + 2p\gamma - \frac{2}{3}\pi^2 p^2 + 4p^2 H_p^{(2)} \\ &\quad + \frac{2(2p^3 - 5p^2 + 7p - 2)}{(p-1)^2} + \mathcal{O}\left(\frac{1}{n^{1-\varepsilon}}\right). \end{aligned}$$

Of course, the case $p = 1$ is trivial, since we have $\mathbb{P}\{Y_{n,1} = 1\} = 1$ due to the fact that the spanning tree of a single node is the node itself.

The case $p = 2$ is of particular interest, since $Y_{n,2}$ is as described earlier, up to an additive constant, the distance $\Delta_n$ between two randomly selected nodes in a binary search tree of size $n$, viz. $\Delta_n \stackrel{\mathcal{L}}{=} Y_{n,2} - 1$. This parameter was studied already in [7], where the asymptotic normality of the distribution was shown by means of a refined contraction method.

As an insightful referee remarks, one could also obtain the Gaussian limit law for $Y_{n,p}$ (without the precision of the order of convergence obtained here) by studying the difference between $X_{n,p}$ and $Y_{n,p}$, which is the length of the path from the root of the binary search tree to the root of the minimal spanning tree. This quantity is very short, for example, it can be shown, that it is zero with probability $1-2/(p+1)$ asymptotically for $n \to \infty$ and $p \geq 2$. Since we gave already a detailed analysis of $Y_{n,p}$ in this section, we will only describe, very briefly, how one could proceed alternatively. It follows by comparing Theorem 4 and Lemma 1, that $\mathbb{E}(X_{n,p} - Y_{n,p}) = 4 + 2p/(p-1) + \mathcal{O}(1/n^{1-\varepsilon})$. One gets thus, that $\mathbb{P}\{X_{n,p} - Y_{n,p} \geq (\log n)^{1/4}\} = \mathcal{O}((\log n)^{-1/4})$. This bound finally suffices to transfer the limiting distribution result from $X_{n,p}$ to $Y_{n,p}$ by considering $\mathbb{P}\{(Y_{n,p} - 2p\log n)/\sqrt{2p\log n} < x\} = \mathbb{P}\{(X_{n,p} - 2p\log n)/\sqrt{2p\log n} - (X_{n,p} - Y_{n,p})/\sqrt{2p\log n} < x\}$.

Institut für Diskrete Mathematik
und Geometrie
Technische Universität Wien
Wiedner Hauptstrasse 8-10
A-1040 Wien
Austria
e-mail: Alois.Panholzer@tuwien.ac.at

The John Knopfmacher Centre
for Applicable Analysis
and Number Theory
School of Mathematics
University of the Witwatersrand
P.O. Wits
2050 Johannesburg
South Africa
e-mail: helmut@staff.ms.wits.ac.za
url: http://www.wits.ac.za/helmut/index.htm